\documentclass{gtart}

\def\ifplaintex{\expandafter\ifx\csname documentclass\endcsname\relax}

\ifplaintex 
\else
\fi

\def\gtm{{\mathsurround=0pt\it $\cal G\mskip-2mu$eometry \&\ 
$\cal T\!\!$opology $\cal M\mskip-1mu$onographs}}    %  for monographs

\def\gtp{{\mathsurround=0pt\it $\cal G\mskip-2mu$eometry \&\ 
$\cal T\!\!$opology $\cal P\!$ublications}}  % GT publications

\def\recd{{\small Received:\qua\receiveddate\ifx\reviseddate\relax
\else\qquad Revised:\qua\reviseddate\fi\par}} 

\def\volumenumber#1{\def\thevolumenumber{#1}}
\def\volumeyear#1{\def\thevolumeyear{#1}}
\def\volumename#1{\def\thevolumename{#1}}
\def\papernumber#1{\def\thepapernumber{#1}}
\def\pagenumbers#1#2{\def\startpage{#1}\def\finishpage{#2}}
\def\published#1{\def\publishdate{#1}}
\def\received#1{\def\receiveddate{#1}}

\def\accepted#1{\def\accepteddate{#1}}

\def\asciiaddress#1{\def\theasciiaddress{#1}}

\long\def\asciiabstract#1{\long\def\theasciiabstract{#1}}

\let\\\par
\let\thevolumenumber\relax\let\thepapernumber\relax
\let\thevolumeyear\relax\let\startpage\relax
\let\finishpage\relax\let\publishdate\relax\let\receiveddate\relax
\let\reviseddate\relax\let\accepteddate\relax\let\theasciititle\relax
\let\theasciiauthors\relax\let\theasciiaddress\relax
\let\theasciiabstract\relax

\let\theerratum\relax\let\theasciiemail\relax
\let\theshortauthors\relax\let\theshorttitle\relax

\def\startpage{1}\def\finishpage{15}\def\thepapernumber{77}

\volumenumber{2}
\volumename{Proceedings of the Kirbyfest}
\volumeyear{1999}

\long\def\maketitlep{   % start of definition of \maketitlep

\count0=\startpage

\gtm\nl        %   GT mongraphs (top left) 
{\small Volume \thevolumenumber: \thevolumename\nl 
\ifx\theerratum\relax\else Erratum \erratumnumber\nl\fi
Pages \startpage--\finishpage\nl}

\vglue 0.1truein   % top margin

{\parskip=0pt\leftskip 0pt plus 1fil\def\\{\par\smallskip}{\ifplaintex\large
\else\Large\fi\bf\thetitle}\par\medskip}   
\vglue 0.05truein 

{\parskip=0pt\leftskip 0pt plus 1fil\def\\{\par}{\sc\theauthors}
\par\medskip}%
 
\vglue 0.03truein 

{\small\leftskip 25pt\rightskip 25pt{\bf Abstract}\stdspace\theabstract

{\bf AMS Classification}\stdspace\theprimaryclass
\ifx\thesecondaryclass\relax\else; \thesecondaryclass\fi\par
{\bf Keywords}\stdspace \thekeywords\par}\vglue 7pt

}   % end of definition of \maketitlep

\font\phead=cmsl9 scaled 950
\font=cmsl9 scaled 1050
\font\pnum=cmbx10 scaled 913
\font\lnum=cmbx10 
\font\pfoot=cmsl9 scaled 950
\font=cmsl9 scaled 1050
\ifplaintex
\headline{\vbox to 0pt{\vskip -4.5mm\line{\small\phead\ifnum
\count0=\startpage ISSN 1464-8997 (on line)
1464-8989 (printed) \hfill {\pnum\folio}\else\ifodd\count0\def\\{ }% 
\ifx\theshorttitle\relax\thetitle\else\theshorttitle\fi\hfill{\pnum\folio}
\else\def\\{ and }{\pnum\folio}\hfill\ifx\theshortauthors\relax\theauthors
\else\theshortauthors\fi\fi\fi}\vss}}
\footline{\vbox to 0pt{\vglue 0mm\line{\small\pfoot\ifnum\count0=\startpage
Published \publishdate:\qua\copyright\ \gtp\hfill\else
\gtm, Volume \thevolumenumber\ (\thevolumeyear)\hfill\fi}\vss
}}
\else
\makeatletter
\def\@oddhead{{\small\lhead\ifnum\count0=\startpage ISSN 1464-8997 (on line)
1464-8989 (printed) \hfill {\lnum\number\count0}\else\ifodd\count0
\def\\{ }\ifx\theshorttitle\relax \thetitle \else\theshorttitle\fi\hfill
{\lnum\number\count0}\else\def\\{ and }{\lnum\number\count0}
\hfill\ifx\theshortauthors\relax 
\theauthors\else\theshortauthors\fi\fi\fi}}\def\@evenhead{@oddhead}
\def\@oddfoot{\small\lfoot\ifnum\count0=\startpage Published \publishdate:\qua\copyright\ \gtp\hfill\else
\gtm, Volume \thevolumenumber\ (\thevolumeyear)\hfill\fi}
\def\@evenfoot{@oddfoot}
\makeatother
\fi

\let\maketitlepage\maketitlep
\let\makeshorttitle\maketitlepage
\let\maketitle\maketitlepage

\newwrite\gtoutfile
\long\gdef\makeheadfile{  %%% start of definition of \makeheadfile
{\def\\{, }\def\s{ }
\immediate\openout\gtoutfile head.xxx
\immediate\write\gtoutfile{To: math@arxiv.org}
\immediate\write\gtoutfile{Subject: put OR rep NNNNN:ppppp}
\immediate\write\gtoutfile{--text follows this line--}
\immediate\write\gtoutfile{Proxy-for: \ifx\theasciiauthors\relax
\theauthors\else\theasciiauthors\fi\s<\ifx\theasciiemail\relax\theemail\else\theasciiemail\fi>}
\immediate\write\gtoutfile{\noexpand\\}
\immediate\write\gtoutfile{Authors: \ifx\theasciiauthors\relax
\theauthors\else\theasciiauthors\fi}
{\def\\{ }\immediate\write\gtoutfile{Title: \ifx\theasciititle\relax
\thetitle\else\theasciititle\fi}}
\immediate\write\gtoutfile{Subj-class: GT or SG, GR etc}
\immediate\write\gtoutfile{MSC-class: \theprimaryclass\ifx\thesecondaryclass\relax\else, \thesecondaryclass\fi}
\immediate\write\gtoutfile{Journal-ref: Geom. Topol. Monogr. \thevolumenumber\s
(\thevolumeyear) \startpage-\finishpage}
\immediate\write\gtoutfile{Comments: Published by Geometry and Topology Monographs at}
\immediate\write\gtoutfile{\s\s\s  http://www.maths.warwick.ac.uk/gt/GTMon\thevolumenumber/paper\thepapernumber.abs.html}
\immediate\write\gtoutfile{\noexpand\\}
\immediate\write\gtoutfile{}
\ifx\theasciiabstract\relax
\immediate\write\gtoutfile{\theabstract}\else
\immediate\write\gtoutfile{\theasciiabstract}\fi
\immediate\write\gtoutfile{}
\immediate\write\gtoutfile{\noexpand\\}
\immediate\write\gtoutfile{}
\immediate\closeout\gtoutfile}}  %%% end of definition of \makeheadfile

\def\maketitlepage{\maketitlep\makeheadfile}
\let\makeshorttitle\maketitlepage
\let\maketitle\maketitlepage

\volumenumber{4}
\volumename{Invariants of knots and 3-manifolds (Kyoto 2001)}
\volumeyear{2002}
\papernumber{13}
\pagenumbers{201}{214}
\received{12 December 2001}
\accepted{22 July 2002}
\published{21 September 2002}

\usepackage{amsmath,amssymb}
\usepackage{pstricks}

\def\lmn#1{\vadjust{\setbox1=\vtop{\hsize 12mm
\parindent=0pt\baselineskip=9pt
\rightskip=4mm plus 4mm#1}
\hbox{\kern-12mm\smash{\raise .5ex\box1}}}}

\newtheorem{theorem}{Theorem}[section]

\newtheorem{proposition}[theorem]{Proposition}

\def\P{\mathcal P}
\def\D{\mathcal D}

\def\s{\sigma}

\def\Y{\mathrm Y}
\def\LP{{\Lambda^{(p)}}}

\def\Z{{\mathbb Z}}

\def\build#1_#2^#3{\mathrel{\mathop{\kern 0pt#1}\limits_{#2}^{#3}}}

\begin{document}

%% GM: \input{figdef5}

\newcommand{\MILN}{
%\begin{pspicture}[.4](-.1,-.1)(1.3,1.1)
\begin{pspicture}[0](0,.1)(.8,.6)
%\psset{unit=.5cm}
%\psgrid
\psset{unit=.5}
\psline[linewidth=.5pt,linestyle=dotted,dotsep=1pt](.1,.1)(.6,.5)
\psline[linewidth=.5pt,linestyle=dotted,dotsep=1pt](.1,.9)(.6,.5)
\psline[linewidth=.5pt,linestyle=dotted,dotsep=1pt](.6,.5)(1.1,.5)
\end{pspicture}
}

\newcommand{\FIGEXi}{
\begin{pspicture}[.4](-.5,1)(4.5,3.5)
%\psgrid
%\psset{unit=1cm}
\psline[linewidth=.5pt,linestyle=dotted,dotsep=1pt](2,1.5)(3,2)
\psline[linewidth=.5pt,linestyle=dotted,dotsep=1pt](2,1.5)(1,2)
\psline[linewidth=.5pt,linestyle=dotted,dotsep=1pt](3,2)(3.5,3)
\psline[linewidth=.5pt,linestyle=dotted,dotsep=1pt](3,2)(4,1.5)
\psline[linewidth=.5pt,linestyle=dotted,dotsep=1pt](1,2)(.5,3)
\psline[linewidth=.5pt,linestyle=dotted,dotsep=1pt](1,2)(0,1.5)
\put(-.4,1.3){$\scriptstyle{3}$}
\put(.1,2.9){$\scriptstyle{2}$}
\put(1.8,.9){$\scriptstyle{1}$}
\put(4.1,1.3){$\scriptstyle{4}$}
\put(3.6,2.9){$\scriptstyle{5}$}
\psdots[dotscale=.6](0,1.5)(.5,3)(2,1.5)(3.5,3)(4,1.5)
\end{pspicture}
}

\newcommand{\MILNijk}{
%\begin{pspicture}[.4](-.1,-.1)(1.3,1.1)
\begin{pspicture}[.4](-.3,-.3)(1.5,1.3)
%\psgrid
\psline[linewidth=.5pt,linestyle=dotted,dotsep=1pt](.1,.1)(.6,.5)
\psline[linewidth=.5pt,linestyle=dotted,dotsep=1pt](.1,.9)(.6,.5)
\psline[linewidth=.5pt,linestyle=dotted,dotsep=1pt](.6,.5)(1.1,.5)
\put(1.3,.4){$\scriptstyle{k}$}
\put(-.2,0){$\scriptstyle{j}$}
\put(-.2,.8){$\scriptstyle{i}$}
\end{pspicture}
}

\def\KPlus{
\begin{picture}(2,2)(-1,-1)
%\put(0,0){\circle{2}}
\put(-0.707,-0.707){\vector(1,1){1.414}}
\put(0.707,-0.707){\line(-1,1){0.6}}
\put(-0.107,0.107){\vector(-1,1){0.6}}
\end{picture}}

\def\KMinus{
\begin{picture}(2,2)(-1,-1)
%\put(0,0){\circle{2}}
\put(-0.707,-0.707){\line(1,1){0.6}}
\put(0.107,0.107){\vector(1,1){0.6}}
\put(0.707,-0.707){\vector(-1,1){1.414}}
\end{picture}}

\def\KII{
\begin{picture}(2,2)(-1,-1)
%\put(0,0){\circle{2}}
\qbezier(-0.707,-0.707)(0,0)(-0.707,0.707)
\qbezier(0.707,-0.707)(0,0)(0.707,0.707)
\put(-0.607,0.607){\vector(-1,1){0.1414}}
\put(0.607,0.607){\vector(1,1){0.1414}}
\end{picture}}

\newcommand{\Wi}{
\begin{pspicture}[.4](0,0)(1,1)
\psline{->}(0,0)(0,1)
\psline{->}(1,1)(1,0)
\psline[linewidth=.5pt,linestyle=dotted,dotsep=1pt](0,.5)(1,.5)
\end{pspicture}
}

\newcommand{\Wii}{
\begin{pspicture}[.4](0,0)(1,1)
\psline{->}(0,0)(0,.4)(1,.4)(1,0)
\psline{->}(1,1)(1,.6)(0,.6)(0,1)
\end{pspicture}
}

\newcommand{\Gi}{
\begin{pspicture}[-.5](0,0)(1,.2)
\psline(0,0)(1,0)
\psdots(0,0)(1,0)
\put(1,-.3){$\scriptstyle{2}$}
\put(0,-.3){$\scriptstyle{1}$}
\end{pspicture}
}

\newcommand{\Gii}{
\begin{pspicture}[.4](0,-.2)(1,1)
\psline(0,0)(1,0)(.5,.85)(0,0)
\psdots(0,0)(1,0)(.5,.85)
\put(1.1,-.3){$\scriptstyle{2}$}
\put(-.2,-.3){$\scriptstyle{1}$}
\put(.45,1){$\scriptstyle{3}$}
\end{pspicture}
}
\newcommand{\Gia}{
\begin{pspicture}[.4](0,-.2)(1,1)
\psline(0,0)(1,0)(.5,.85)
\psdots(0,0)(1,0)(.5,.85)
\put(1.1,-.3){$\scriptstyle{2}$}
\put(-.2,-.3){$\scriptstyle{1}$}
\put(.45,1){$\scriptstyle{3}$}
\end{pspicture}
}
\newcommand{\Gib}{
\begin{pspicture}[.4](0,-.2)(1,1)
\psline(1,0)(.5,.85)(0,0)
\psdots(0,0)(1,0)(.5,.85)
\put(1.1,-.3){$\scriptstyle{2}$}
\put(-.2,-.3){$\scriptstyle{1}$}
\put(.45,1){$\scriptstyle{3}$}
\end{pspicture}
}
\newcommand{\Gic}{
\begin{pspicture}[.4](0,-.2)(1,1)
\psline(.5,.85)(0,0)(1,0)
\psdots(0,0)(1,0)(.5,.85)
\put(1.1,-.3){$\scriptstyle{2}$}
\put(-.2,-.3){$\scriptstyle{1}$}
\put(.45,1){$\scriptstyle{3}$}
\end{pspicture}
}

\newcommand{\ACE}{
\begin{pspicture}[.5](0,-.2)(3,1)
\pscircle(.5,.5){.35}
\pscircle(1.5,.5){.35}
\pscircle(2.5,.5){.35}
\psarc[arrowlength=.7,arrowsize=2pt 4]{->}(.5,.5){.35}{20}{95}
\psarc[arrowlength=.7,arrowsize=2pt 4]{->}(1.5,.5){.35}{20}{95}
\psarc[arrowlength=.7,arrowsize=2pt 4]{->}(2.5,.5){.35}{20}{95}
\psline[linewidth=.5pt,linestyle=dotted,dotsep=1pt](.85,.5)(1.15,.5)
\psline[linewidth=.5pt,linestyle=dotted,dotsep=1pt](1.85,.5)(2.15,.5)
\put(.85,.2){$\scriptstyle{1}$}
\put(1.85,.2){$\scriptstyle{2}$}
\put(2.85,.2){$\scriptstyle{3}$}
\end{pspicture}
}

\newcommand{\ACEi}{
\begin{pspicture}[.5](0,-.2)(3,1)
%\psgrid
%\psset{unit=4cm}
\psarc(.5,.5){.35}{20}{340}
\psarc(1.5,.5){.35}{20}{160}
\psarc(1.5,.5){.35}{200}{340}
\psarc(2.5,.5){.35}{200}{160}
\psarc[arrowlength=.7,arrowsize=2pt 4]{->}(1.5,.5){.35}{20}{95}
\psline(.825,.385)(1.175,.385)
\psline(1.825,.385)(2.175,.385)
\psline(.825,.615)(1.175,.615)
\psline(1.825,.615)(2.175,.615)
\end{pspicture}
}
\newcommand{\Tree}{
\begin{pspicture}[-.6](0,0)(1,.2)
\psline(0,0)(1,0)
\psdots(0,0)(.5,0)(1,0)
\put(0,-.3){$\scriptstyle{1}$}
\put(.5,-.3){$\scriptstyle{2}$}
\put(1,-.3){$\scriptstyle{3}$}
\end{pspicture}
}

\newcommand{\Stepa}{
\begin{pspicture}[.4](0,0)(1,1)
\pscircle(.5,.5){.35}
\psarc[arrowlength=.7,arrowsize=2pt 4]{->}(.5,.5){.35}{20}{95}
\psarc{->}(1.5,.5){.35}{120}{240}
\psline[linewidth=.5pt,linestyle=dotted,dotsep=1pt](.85,.5)(1.15,.5)
\end{pspicture}
}
\newcommand{\Stepai}{
\begin{pspicture}[.4](0,0)(1,1)
\psarc{->}(.5,.5){.35}{120}{240}
\end{pspicture}
}

\newcommand{\Stepaa}{
\begin{pspicture}[.4](0,0)(1,1)
%\psset{unit=4cm}
%\psgrid
\pscircle(.5,.5){.35}
\psarc[arrowlength=.7,arrowsize=2pt 4]{->}(.5,.5){.35}{20}{95}
\psarc{->}(1.7,.5){.35}{120}{240}
\psline[linewidth=.5pt,linestyle=dotted,dotsep=1pt](.85,.5)(1.05,.5)(1.4,.3)
\psline[linewidth=.5pt,linestyle=dotted,dotsep=1pt](1.05,.5)(1.4,.7)
\end{pspicture}
}

\newcommand{\Stepb}{
\begin{pspicture}[.4](0,0)(1.5,1)
%\psset{unit=4cm}
%\psgrid
\pscircle(.5,.5){.35}
\psarc[arrowlength=.7,arrowsize=2pt 4]{->}(.5,.5){.35}{20}{95}
\psarc{->}(1.5,.25){.35}{150}{210}
\psarc{->}(1.5,.75){.35}{150}{210}
\psline[linewidth=.5pt,linestyle=dotted,dotsep=1pt](.8,.68)(1.15,.74)
\psline[linewidth=.5pt,linestyle=dotted,dotsep=1pt](.8,.32)(1.15,.26)
\end{pspicture}
}

\newcommand{\Stepbi}{
\begin{pspicture}[.4](0,0)(1.5,1)
%\psset{unit=4cm}
%\psgrid
\psarc{->}(1.5,.25){.35}{150}{210}
\psarc{->}(1.5,.75){.35}{150}{210}
\pscurve[linewidth=.5pt,linestyle=dotted,dotsep=1pt](1.15,.74)(.9,.6)(.9,.4)(1.15,.26)
\end{pspicture}
}

\newcommand{\Stepc}{
\begin{pspicture}[.4](0,0)(1.5,1)
%\psset{unit=4cm}
%\psgrid
\pscircle(.5,.5){.35}
\psarc[arrowlength=.7,arrowsize=2pt 4]{->}(.5,.5){.35}{20}{95}
\psarc{->}(1.8,.25){.35}{150}{210}
\psarc{->}(1.8,.75){.35}{150}{210}
\psline[linewidth=.5pt,linestyle=dotted,dotsep=1pt](.85,.5)(1.15,.5)(1.45,.26)
\psline[linewidth=.5pt,linestyle=dotted,dotsep=1pt](1.15,.5)(1.45,.74)
\psarc{->}(-.7,.5){.35}{300}{60}
\psline[linewidth=.5pt,linestyle=dotted,dotsep=1pt](-.35,.5)(.15,.5)
\end{pspicture}
}
\newcommand{\Stepci}{
\begin{pspicture}[.4](0,0)(1.5,1)
\psarc{->}(1.8,.25){.35}{150}{210}
\psarc{->}(1.8,.75){.35}{150}{210}
\psline[linewidth=.5pt,linestyle=dotted,dotsep=1pt](.85,.5)(1.15,.5)(1.45,.26)
\psline[linewidth=.5pt,linestyle=dotted,dotsep=1pt](1.15,.5)(1.45,.74)
\psarc{->}(.3,.5){.35}{300}{60}
\psline[linewidth=.5pt,linestyle=dotted,dotsep=1pt](.65,.5)(1.15,.5)
\end{pspicture}
}

\newcommand{\Stepd}{
\begin{pspicture}[.4](0,0)(1.5,1)
%\psset{unit=4cm}
%\psgrid
\pscircle(.5,.5){.35}
\psarc[arrowlength=.7,arrowsize=2pt 4]{->}(.5,.5){.35}{20}{95}
\psarc{->}(1.8,.25){.35}{150}{210}
\psarc{->}(1.8,.75){.35}{150}{210}
\psline[linewidth=.5pt,linestyle=dotted,dotsep=1pt](.85,.5)(1.15,.5)(1.45,.26)
\psline[linewidth=.5pt,linestyle=dotted,dotsep=1pt](1.15,.5)(1.45,.74)
\psarc{->}(-.8,.25){.35}{330}{30}
\psarc{->}(-.8,.75){.35}{330}{30}
\psline[linewidth=.5pt,linestyle=dotted,dotsep=1pt](.15,.5)(-.15,.5)(-.45,.26)
\psline[linewidth=.5pt,linestyle=dotted,dotsep=1pt](-.15,.5)(-.45,.74)
\end{pspicture}
}

\newcommand{\Stepdi}{
\begin{pspicture}[.4](0,0)(1.5,1)
\psarc{->}(1.8,.25){.35}{150}{210}
\psarc{->}(1.8,.75){.35}{150}{210}
\psline[linewidth=.5pt,linestyle=dotted,dotsep=1pt](.85,.5)(1.15,.5)(1.45,.26)
\psline[linewidth=.5pt,linestyle=dotted,dotsep=1pt](1.15,.5)(1.45,.74)
\psarc{->}(.2,.25){.35}{330}{30}
\psarc{->}(.2,.75){.35}{330}{30}
\psline[linewidth=.5pt,linestyle=dotted,dotsep=1pt](.85,.5)(.55,.26)
\psline[linewidth=.5pt,linestyle=dotted,dotsep=1pt](.85,.5)(.55,.74)
\end{pspicture}
}

\newcommand{\Stepdii}{
\begin{pspicture}[.4](0,0)(1.5,1)
\psarc{->}(1.8,.25){.35}{150}{210}
\psarc{->}(1.8,.75){.35}{150}{210}
\psline[linewidth=.5pt,linestyle=dotted,dotsep=1pt](.55,.74)(1.45,.74)
\psarc{->}(.2,.25){.35}{330}{30}
\psarc{->}(.2,.75){.35}{330}{30}
\end{pspicture}
}

\newcommand{\Stepdiia}{
\begin{pspicture}[.4](0,0)(1.5,1)
\psarc{->}(1.8,.25){.35}{150}{210}
\psarc{->}(1.8,.75){.35}{150}{210}
\psline[linewidth=.5pt,linestyle=dotted,dotsep=1pt](.55,.26)(1.45,.26)
\psarc{->}(.2,.25){.35}{330}{30}
\psarc{->}(.2,.75){.35}{330}{30}
\end{pspicture}
}

\newcommand{\Stepdiib}{
\begin{pspicture}[.4](0,0)(1.5,1)
\psarc{->}(1.8,.25){.35}{150}{210}
\psarc{->}(1.8,.75){.35}{150}{210}
\psline[linewidth=.5pt,linestyle=dotted,dotsep=1pt](.55,.74)(1.45,.26)
\psarc{->}(.2,.25){.35}{330}{30}
\psarc{->}(.2,.75){.35}{330}{30}
\end{pspicture}
}

\newcommand{\Stepdiic}{
\begin{pspicture}[.4](0,0)(1.5,1)
\psarc{->}(1.8,.25){.35}{150}{210}
\psarc{->}(1.8,.75){.35}{150}{210}
\psline[linewidth=.5pt,linestyle=dotted,dotsep=1pt](.55,.26)(1.45,.74)
\psarc{->}(.2,.25){.35}{330}{30}
\psarc{->}(.2,.75){.35}{330}{30}
\end{pspicture}
}

\newcommand{\ExThree}{
\begin{pspicture}[.3](-.1,-.6)(2,1.1)
%\psset{unit=4cm}
%\psgrid
\pscircle(0,.5){.35}
\psarc[arrowlength=.7,arrowsize=2pt 4]{->}(0,.5){.35}{20}{95}
\put(-.1,.56){$\scriptstyle 1$}
\pscircle(1.5,.5){.35}
\put(1.4,.56){$\scriptstyle 3$}
\psarc[arrowlength=.7,arrowsize=2pt 4]{->}(1.5,.5){.35}{20}{95}
\pscircle(.75,-.4){.35}
\put(.65,-.34){$\scriptstyle 2$}
\psarc[arrowlength=.7,arrowsize=2pt 4]{->}(.75,-.4){.35}{200}{275}
\psline[linewidth=.5pt,linestyle=dotted,dotsep=1pt](.3,.7)(1.2,.7)
\psline[linewidth=.5pt,linestyle=dotted,dotsep=1pt](.3,.3)(1.2,.3)
\psline[linewidth=.5pt,linestyle=dotted,dotsep=1pt](.6,.7)(.6,-.1)
\psline[linewidth=.5pt,linestyle=dotted,dotsep=1pt](.9,.3)(.9,-.1)
\end{pspicture}
}
\newcommand{\ExThreei}{
\begin{pspicture}[.6](-.5,-.6)(2,1.1)
%%  \psset{unit=4cm}
%% \psgrid
\pscircle(0,.5){.35}
\put(-.1,.56){$\scriptstyle 1$}
\psarc[arrowlength=.7,arrowsize=2pt 4]{->}(0,.5){.35}{20}{95}
\pscircle(1.5,.5){.35}
\put(1.4,.56){$\scriptstyle 3$}
\psarc[arrowlength=.7,arrowsize=2pt 4]{->}(1.5,.5){.35}{20}{95}
\psline[linewidth=.5pt,linestyle=dotted,dotsep=1pt](.35,.5)(1.15,.5)
\end{pspicture}
}
\newcommand{\ExThreeii}{
\begin{pspicture}[.6](-.5,-.6)(.9,1.1)
%\psset{unit=4cm}
%\psgrid
\pscircle(0,.5){.35}
\put(-.1,.56){$\scriptstyle 1$}
\psarc[arrowlength=.7,arrowsize=2pt 4]{->}(0,.5){.35}{20}{95}
\end{pspicture}
}

\newcommand{\ExFive}{
\begin{pspicture}[.6](-.5,-1.5)(2,1)
%\psset{unit=4cm}
%\psgrid
\pscircle(0,.5){.35}
\put(-.1,.56){$\scriptstyle 1$}
\psarc[arrowlength=.7,arrowsize=2pt 4]{->}(0,.5){.35}{20}{95}
\pscircle(1.5,.5){.35}
\put(1.4,.56){$\scriptstyle 5$}
\psarc[arrowlength=.7,arrowsize=2pt 4]{->}(1.5,.5){.35}{20}{95}
\pscircle(-.5,-.5){.35}
\put(-.6,-.44){$\scriptstyle 2$}
\psarc[arrowlength=.7,arrowsize=2pt 4]{->}(-.5,-.5){.35}{20}{95}
\pscircle(2,-.5){.35}
\put(1.9,-.44){$\scriptstyle 4$}
\psarc[arrowlength=.7,arrowsize=2pt 4]{->}(2,-.5){.35}{20}{95}
\pscircle(.75,-1.2){.35}
\put(.65,-1.14){$\scriptstyle 3$}
\psarc[arrowlength=.7,arrowsize=2pt 4]{->}(.75,-1.2){.35}{20}{95}
\pscurve[linewidth=.5pt,linestyle=dotted,dotsep=1pt](.35,.5)(.75,.4)(1.15,.5)
\pscurve[linewidth=.5pt,linestyle=dotted,dotsep=1pt](1.7,-.3)(1,0)(.75,.4)
\pscurve[linewidth=.5pt,linestyle=dotted,dotsep=1pt](.2,.2)(.1,-.2)(-.2,-.35)
\psline[linewidth=.5pt,linestyle=dotted,dotsep=1pt](.46,-1)(.1,-.2)
\pscurve[linewidth=.5pt,linestyle=dotted,dotsep=1pt](1.6,0.17)(1.63,0)(1.8,-0.2)
\pscurve[linewidth=.5pt,linestyle=dotted,dotsep=1pt](1.63,0)(1,-.5)(.9,-.9)
\pscurve[linewidth=.5pt,linestyle=dotted,dotsep=1pt](1.3,.2)(.6,-.4)(-.16,-.6)
\psline[linewidth=.5pt,linestyle=dotted,dotsep=1pt](.6,-.4)(.7,-.85)
\end{pspicture}
}

\newcommand{\ExFour}{
\begin{pspicture}[.6](-.5,0)(2,1)
%%\psset{unit=4cm}
%\psgrid
\pscircle(0,.5){.35}
\put(-.1,.56){$\scriptstyle 1$}
\psarc[arrowlength=.7,arrowsize=2pt 4]{->}(0,.5){.35}{20}{95}
\pscircle(1.5,.5){.35}
\put(1.4,.56){$\scriptstyle 2$}
\psarc[arrowlength=.7,arrowsize=2pt 4]{->}(1.5,.5){.35}{20}{95}
\psline[linewidth=.5pt,linestyle=dotted,dotsep=1pt](.3,.7)(1.2,.7)
\psline[linewidth=.5pt,linestyle=dotted,dotsep=1pt](.3,.3)(1.2,.3)
\psline[linewidth=.5pt,linestyle=dotted,dotsep=1pt](.75,.7)(.75,.3)
\end{pspicture}
}

\newcommand{\ExFouri}{
\begin{pspicture}[.6](-.5,-1)(2,1)
%\psset{unit=4cm}
%\psgrid
\pscircle(0,.5){.35}
\put(-.1,.56){$\scriptstyle 1$}
\psarc[arrowlength=.7,arrowsize=2pt 4]{->}(0,.5){.35}{20}{95}
\pscircle(1.5,.5){.35}
\put(1.4,.56){$\scriptstyle 4$}
\psarc[arrowlength=.7,arrowsize=2pt 4]{->}(1.5,.5){.35}{20}{95}
\psline[linewidth=.5pt,linestyle=dotted,dotsep=1pt](.3,.7)(1.2,.7)
\psline[linewidth=.5pt,linestyle=dotted,dotsep=1pt](.3,.3)(1.2,.3)
\psline[linewidth=.5pt,linestyle=dotted,dotsep=1pt](.75,.7)(.75,.3)
\pscircle(0,-.8){.35}
\put(-.1,-.74){$\scriptstyle 2$}
\psarc[arrowlength=.7,arrowsize=2pt 4]{->}(0,-.8){.35}{200}{275}
\pscircle(1.5,-.8){.35}
\put(1.4,-.74){$\scriptstyle 3$}
\psarc[arrowlength=.7,arrowsize=2pt 4]{->}(1.5,-.8){.35}{200}{275}
\psline[linewidth=.5pt,linestyle=dotted,dotsep=1pt](1.5,.15)(1.5,-.45)
\psline[linewidth=.5pt,linestyle=dotted,dotsep=1pt](0,.15)(0,-.45)
\pscurve[linewidth=.5pt,linestyle=dotted,dotsep=1pt](0,-.1)(1,-.3)(1.27,-.53)
\pscurve[linewidth=.5pt,linestyle=dotted,dotsep=1pt](1.5,-.1)(.5,-.3)(.23,-.53)
\end{pspicture}
}

\newcommand{\Pf}{\operatorname{Pf}}

\title{Matrix-tree theorems and  the\\Alexander-Conway polynomial}
\author{Gregor Masbaum}

\address{Institut de Math{\'e}matiques de Jussieu, Universit{\'e} 
Paris VII\\Case 7012, 75251 Paris Cedex 05, France}
\asciiaddress{Institut de Mathematiques de Jussieu, Universite 
Paris VII\\Case 7012, 75251 Paris Cedex 05, France}
\email{masbaum@math.jussieu.fr}

\begin{abstract} 
This talk is a report on joint work with A.~Vaintrob
\cite{MV1,MV2}. It is organised as follows. We begin by recalling how
the classical Matrix-Tree Theorem relates two different expressions
for the lowest degree coefficient of the Alexander-Conway polynomial
of a link. We then state our formula for the lowest degree coefficient
of an algebraically split link in terms of Milnor's triple linking
numbers.  We explain how this formula can be deduced from a
determinantal expression due to Traldi and Levine by means of our
Pfaffian Matrix-Tree Theorem \cite{MV1}. We also discuss the approach
via finite type invariants, which allowed us in \cite{MV2} to obtain
the same result directly from some properties of the Alexander-Conway
weight system. This approach also gives similar results if all Milnor
numbers up to a given order vanish.
\end{abstract}

\asciiabstract{ 
This talk is a report on joint work with A. Vaintrob
[arXiv:math.CO/0109104 and math.GT/0111102]. It is organised as
follows. We begin by recalling how the classical Matrix-Tree Theorem
relates two different expressions for the lowest degree coefficient of
the Alexander-Conway polynomial of a link. We then state our formula
for the lowest degree coefficient of an algebraically split link in
terms of Milnor's triple linking numbers.  We explain how this formula
can be deduced from a determinantal expression due to Traldi and
Levine by means of our Pfaffian Matrix-Tree Theorem
[arXiv:math.CO/0109104]. We also discuss the approach via finite type
invariants, which allowed us in [arXiv:math.GT/0111102] to obtain the
same result directly from some properties of the Alexander-Conway
weight system. This approach also gives similar results if all Milnor
numbers up to a given order vanish.}


%% GM: \date{}
%\date{October 2001}

%\maketitle
\maketitlepage

\section{The Alexander-Conway polynomial and its lowest order coefficient}
\begin{figure}[ht!]
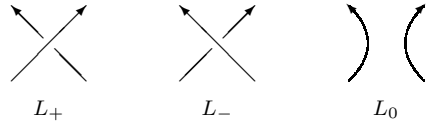

\begin{center}
\begin{displaymath}
\mathop{\KPlus}_{L_+}\qquad
\mathop{\KMinus}_{L_-}\qquad
\mathop{\KII}_{L_0}
\end{displaymath}
\vspace{-20pt}
\caption{\label{skeintriple} A skein triple}
\end{center}
\end{figure}

Let $L$ be an
oriented link in  $S^3$ with $m$ (numbered) components. Its Alexander-Conway polynomial 
$$\nabla_L(z)= \sum_{i\geq 0} c_i(L) z^i \in \Z[z]$$ 
is one of the most thoroughly studied classical isotopy invariants of 
links. It can be defined in various ways. For example, if $V$ is a Seifert matrix for $L$, then 
\begin{equation}
  \label{eq:1}
  \nabla_L(z)= \det (t V - t^{-1} V^T)
\end{equation}
where $z=t-t^{-1}$. Another definition is via the {\em skein
  relation} 
\begin{equation}
\nabla_{L_+}-\nabla_{L_-}=z \nabla_{L_0}~,\label{skein}
\end{equation} 
where $(L_+, L_-, L_0)$ is any skein triple (see
Figure~\ref{skeintriple}).  

Indeed, the Alexander-Conway polynomial is uniquely determined by
the skein relation~(\ref{skein}) and the initial conditions
\begin{equation}
\nabla_{U_m}=\begin{cases} 1 \ \ \ \text{if $m=1$}\\
0\ \ \ \text{if $m\geq 2$},
\end{cases}
\label{unlink}
\end{equation}
where $U_m$ is the trivial link with $m$ components.

Hosokawa~\cite{Hw}, Hartley~\cite[(4.7)]{Ha},  
and Hoste~\cite{Ho} showed that the coefficients $c_i(L)$ of
$\nabla_L$ for an $m$-component link $L$ vanish when $i\leq m-2$ 
and that the coefficient $c_{m-1}(L)$ depends only
on the linking numbers $\ell_{ij}(L)$ between the $i$th and $j$th 
components of $L$. Namely,
\begin{equation}
  \label{eq:ho_det}
c_{m-1}(L)= \det\Lambda^{(p)},
\end{equation}
where 
$\Lambda=(\lambda_{ij})$ is the matrix formed by linking numbers
\begin{equation}\label{eq:Lambda}
\lambda_{ij}=
\begin{cases}
\ \ \ \ \ -\ell_{ij}(L), & \text{\ if \ $i\ne j$\ }  \\
\sum_{k\neq i} \ell_{ik}(L), & \text{\ if \ $i=j$}
\end{cases}
\end{equation}
and $\Lambda^{(p)}$ denotes the matrix obtained by removing from
$\Lambda$ the $p$th row and column (it is easy to see that \ 
$\det\Lambda^{(p)}$ \ does not depend on $p$).

Formula \eqref{eq:ho_det} can be proved using the Seifert matrix definition \eqref{eq:1} of $\nabla_L$. We will not give the proof here, but let us indicate how linking numbers come in from this point of view.  Let $\Sigma$ be a Seifert surface for $L$. The key point is  that the Seifert form restricted to $H_1(\partial \Sigma;\Z) \subset H_1(\Sigma;\Z)$ is just given by the linking numbers $\ell_{ij}$.  In particular, for an appropriate choice of basis for $H_1(\Sigma;\Z)$, the Seifert matrix $V$ contains the matrix $\Lambda^{(p)}$ as a submatrix, which then leads to Formula \eqref{eq:ho_det}.

Hartley and Hoste also gave a second expression for
$c_{m-1}(L)$ as a sum over trees: 
\begin{equation}\label{eq:ho_tree}
c_{m-1}(L)= \sum_{T} \prod_{\{i,j\}\in edges(T)} \ell_{ij}(L)~,
\end{equation}
 where $T$ runs through the spanning trees 
in the complete graph $K_m$. (The complete graph $K_m$ has vertices $\{1,2,\ldots,m\}$, and one and only one edge for every unordered pair $\{i,j\}$ of distinct vertices.)

\begin{figure}[ht!]
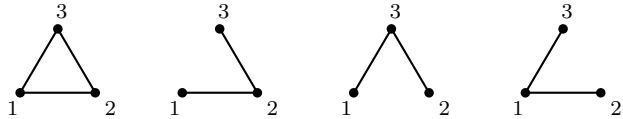

\begin{center}
\Gii \ \ \ \ \ \ \ \Gia\ \ \ \ \ \ \ \ \Gib\ \ \ \ \ \ \ \ \Gic
\caption{\label{Gamma3}The complete graph $K_3$ and its three spanning
trees} 
\end{center}
\end{figure}

For example, if $m=2$ then $c_1(L)=\ell_{12}(L)$, corresponding to 
the only spanning tree in 
 $$K_2=\ \ \ \Gi$$ If $m=3$, then
$$c_2(L)=\ell_{12}(L)\ell_{23}(L) +
\ell_{23}(L)\ell_{13}(L)+\ell_{13}(L)\ell_{12}(L)~,$$ 
corresponding to the three spanning
trees of $K_3$ (see Figure \ref{Gamma3}). 
\vskip 8pt

\section{The classical Matrix-Tree Theorem} \label{sec:2} It is a pleasant exercise to check by hand that Formulas~\eqref{eq:ho_det} 
and~\eqref{eq:ho_tree} give the same answer for $m=2$ and $m=3$. For general $m$ this equality can be deduced from the classical Matrix-Tree Theorem applied to the complete graph $K_m$. 

The statement of this theorem is as follows.  Consider a finite graph $G$ with vertex set $V$ and set of edges
$E$.  If we label each edge $e\in E$  by a variable $x_e$, then a subgraph
of $G$ given as a collection of edges $S\subset E$ corresponds to the
monomial  $$ x_S=\prod_{e\in S}x_e.$$ 
Form a symmetric matrix $\Lambda(G)=(\lambda_{ij})$,  whose rows and columns 
are indexed by the vertices of the graph and entries given by
$$\lambda_{ij}=-\sum_{{e\in E,} \atop {v(e)=\{i,j\} }} x_e, \text{\ \ if\
} i\ne j,   \text{\ \ \ and\ \ \ } 
 \lambda_{ii}=\sum_{{e\in E,} \atop { i\in v(e)}}x_e .$$ Here, we denote by
$v(e) \subset V$ the set of endpoints of the edge $e$.
Since the entries in each row of $\Lambda(G)$ add up to zero, the
determinant of this matrix vanishes and the determinant of the
submatrix $\Lambda(G)^{(p)}$ obtained by deleting the $p$th row and
column of $\Lambda(G)$ is independent of $p$. 
This gives  a polynomial 
\begin{equation}\label{eq:Kir}
\mathcal{D}_G=\det \Lambda(G)^{(p)}
\end{equation} in variables $x_e$ which is
called the Kirchhoff polynomial of $G$. 
The \emph{Matrix-Tree Theorem}~\cite[Theorem VI.29]{Tutte}\cite[Theorem
II.12]{Bol} states that this
polynomial is the generating function of spanning subtrees of the graph $G$ ({\em i.e.} connected acyclic
subgraphs of $G$ with vertex set $V$). 
In other words, one has
\begin{equation}
  \label{eq:mtt}
  \mathcal{D}_G = \sum_T x_T~,
\end{equation}
where the sum is taken over all the spanning subtrees in $G$. 

In the case of the complete graph $K_m$, let us denote its Kirchhoff polynomial by $\mathcal{D}_m$. If we write $x_{ij}=x_{ji}$ for the indeterminate $x_e$ corresponding to the edge $e=\{i,j\}$, then $\Lambda(K_m)$ becomes identified with the matrix $\Lambda$ defined in \eqref{eq:Lambda}. Formula \eqref{eq:ho_det} says that the coefficient $c_{m-1}(L)$ is the Kirchhoff polynomial $\mathcal{D}_m$ evaluated at $x_{ij}=\ell_{ij}(L)$, while Formula \eqref{eq:ho_tree} says that $c_{m-1}(L)$ is the generating function of spanning trees of $K_m$, again  evaluated at $x_{ij}=\ell_{ij}(L)$. Thus, the Matrix-Tree Theorem \eqref{eq:mtt} applied to the complete graph $K_m$ shows that these two formulas for $c_{m-1}(L)$ are equivalent.

\section{An interpretation  via finite type invariants} \label{sec:3}
Formula \eqref{eq:ho_tree} can also be proved directly by induction on the number of components of $L$ \cite{Ha,Ho}. This proof can be  formulated nicely in the language of finite type (Vassiliev) invariants, as follows. (See \cite{MV2} for more details.) 

Recall that the coefficient $c_n$ of the Alexander-Conway
polynomial is a finite type invariant of degree $n$. Let us denote its weight system by $W_n$. It can be computed recursively by  the  
formula 
\begin{equation}
W_{n}(\ \ \Wi\ \ ) = W_{{n-1}}(\ \ \Wii\ \ ) \label{smooth}
\end{equation} 
which follows immediately from the skein relation~(\ref{skein}). For 
a chord diagram $D$ on $m$ circles with $n$ chords, let $D'$ be the result of smoothing 
of all chords by means of~(\ref{smooth}). If $D'$ consists of just one circle, then $$W_{n}(D)=W_{0}(D')=1~.$$ Otherwise, one has $W_{n}(D)=W_{0}(D')=0$.

To see how this relates to 
formula~(\ref{eq:ho_tree}),  note that smoothing of a chord cannot reduce the number of
circles by more than one. Thus, for $W_{n}(D)$ to be non-zero we need at least $m-1$ chords. 
Moreover, the diagrams $D$ with exactly $m-1$ chords 
satisfying $W_{m-1}(D)\ne 0$
must have the property that if each circle
of $D$ is shrinked to a  point, the resulting graph formed by the
chords 
is a tree. See Figure~\ref{Wtree} for an example of 
a chord diagram $D$ whose associated graph is the tree \Tree~.  

\begin{figure}[ht!]
\begin{center}
$W_{2}(\ACE)=W_{0}(\ACEi)=1$
\caption{\label{Wtree}A degree $2$ chord diagram $D$ with $W_{2}(D)=1$}
\end{center}
\end{figure} In other words, the weight system $W_{{m-1}}$ takes the value $1$ on
   precisely those chord 
  diagrams whose associated graph is a spanning tree on the complete
   graph  $K_m$, and $W_{{m-1}}$ is zero on all other chord diagrams.

This simple observation implies Formula~\eqref{eq:ho_tree}, as follows. 
The linking number $\ell_{ij}$ is a finite type invariant of order $1$ 
whose weight system 
is the linear form dual to the chord diagram having just
one chord connecting the $i$th and $j$th circle. It follows that
the right hand side of (\ref{eq:ho_tree}) (which is the spanning tree generating function of $K_m$ evaluated in the $\ell_{ij}$'s) is a finite type invariant of order $m-1$ whose weight system is equal to $W_{{m-1}}$. This proves Formula~(\ref{eq:ho_tree}) on the level of weight systems. The proof can be completed using the fact that the Alexan\-der-Conway polynomial is 
(almost) a 
 canonical invariant \cite{BNG} (see \cite{MV2}).

\section{Algebraically split links and Levine's formula}

If the link $L$ is \emph{algebraically split}, \emph{i.e.}\ all linking
numbers  $\ell_{ij}$ vanish, then not only  $c_{m-1}(L)=0$, but, as
was proved by Traldi~\cite{Tr1,Tr2} and Levine~\cite{Le1}, the next $m-2$ 
coefficients of $\nabla_L$ also vanish
$$c_{m-1}(L)=c_{m}(L)=\ldots=c_{2m-3}(L)=0.$$
For algebraically split oriented links,
there exist well-defined integer-valued isotopy invariants
$\mu_{ijk}(L)$
called the \emph{Milnor triple linking numbers}.
These invariants  generalize  ordinary linking numbers, 
but unlike $\ell_{ij}$, the triple linking numbers  are antisymmetric
with respect to their indices,
$\mu_{ijk}(L)=-\mu_{jik}(L)=\mu_{jki}(L).$  
Thus, for an algebraically split link $L$ with $m$ components,   
we have $m \choose 3$ triple linking numbers $\mu_{ijk}(L)$
corresponding  to the different $3$-component sublinks of $L$.

Levine~\cite{Le1} (see also Traldi~\cite[Theorem~8.2]{Tr2}) found an 
expression for the 
coefficient $c_{2m-2}(L)$ of $\nabla_L$ for an algebraically split 
$m$-component  link in terms  of triple Milnor numbers
\begin{equation}
  \label{eq:lev}
c_{2m-2}(L)=\det\LP,
\end{equation}
where $\Lambda=(\lambda_{ij})$ is an $m\times m$ skew-symmetric matrix
with entries
\begin{equation}
  \label{eq:skew} \lambda_{ij}=\sum_{k} \mu_{ijk}(L),
\end{equation} 
and $\Lambda^{(p)}$, as before, is the result of removing the $p$th 
row and column. 

For example, if $m=3$, 
we have
$$
\Lambda = \left( \begin{array}{ccc}
0 &  \mu_{123}(L)&  \mu_{132}(L)\\
 \mu_{213}(L) & 0&  \mu_{231}(L)\\
 \mu_{312}(L) &  \mu_{321}(L) & 0
\end{array}\right) 
$$
and 
\begin{equation}
  \label{Coc}
  c_4(L)=\det \Lambda^{(3)} = - \mu_{123}(L) \mu_{213}(L)= \mu_{123}(L)^2~.
\end{equation}
 This formula (in the $m=3$ case) goes back to  Cochran~\cite[Theorem 5.1]{Co}.

Similar to Formula~(\ref{eq:ho_det}) for $c_{m-1}(L)$, Levine's proof of Formula~(\ref{eq:lev}) uses the Seifert matrix definition (\ref{eq:1}) of $\nabla_L$. 

\section{The Pfaffian-tree polynomial $\P_m$}

Formula~(\ref{eq:lev}) is 
similar to the first determinantal expression~(\ref{eq:ho_det}). One of the main results of \cite{MV1,MV2} is that there is  an analog of the tree sum formula~(\ref{eq:ho_tree}) for
algebraically split links. 
To state this result, we need to introduce another tree-generating polynomial
analogous to the Kirchhoff polynomial.

Namely, instead of usual graphs whose edges can be thought of as
 segments  joining pairs of points, we consider \emph{$3$-graphs}  
 whose edges have three (distinct) vertices and can be visualized as
 triangles or Y-shaped objects \MILN with the three vertices at their
 endpoints.

\begin{figure}[ht!]
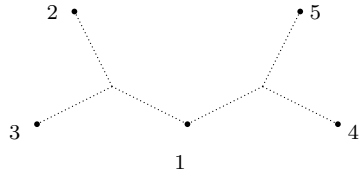

\begin{center}
\FIGEXi
\caption{\label{figexi} A spanning tree in the complete $3$-graph
$\Gamma_5$. It   has two edges, $\{1,2,3\}$ and $\{1,4,5\}$, and
contributes the term  $ y_{123}\, y_{145}$ to $\P_5$.}  
\end{center}
\end{figure}

The notion of spanning trees on a $3$-graph is defined in the natural way.
A sub-$3$-graph $T$ of a $3$-graph $G$ is {\em  spanning} if its vertex set equals that of $G$, and it is a tree if its topological realization (\emph{i.e.}\ the $1$-complex obtained by gluing 
together  Y-shaped objects 
$ \MILN$
corresponding to the edges of $T$) is a tree (\emph{i.e.}\ it is connected
 and  simply connected). See  Figure~\ref{figexi} for an example.

Similarly to the variables $x_{ij}$ of $\D_m$, for each triple of distinct
numbers $i,j,k \in \{1,2,\ldots,m\}$ we introduce   variables
$y_{ijk}$ antisymmetric in $i,j,k$ 
$$
  y_{ijk}=-   y_{jik}=   y_{jki}, \text{\ and\ } y_{iij}=0~.
$$
These variables correspond to edges $\{i,j,k\}$ of the \emph{complete
  $3$-graph} $\Gamma_m$ with vertex set $\{1,\ldots,m\}$.
 
As in the case of ordinary graphs, the correspondence
\begin{equation*}
\text{variable \ }    y_{ijk} \quad  \mapsto \quad  \text{edge\ }
\{i,j,k\} \  \text{of \ } \Gamma_m
\end{equation*}
assigns to each monomial in $  y_{ijk}$  a  sub-$3$-graph
of $\Gamma_m$.

The generating function of spanning trees in the complete $3$-graph $\Gamma_m$
is called the \emph{Pfaffian-tree polynomial} $\P_m$ in \cite{MV1,MV2}. It is 
$$\P_m= \sum_T y_T$$ where the sum is over all spanning trees $T$ of 
$\Gamma_m$, and $y_T$ is, up to sign, just the product of the variables $y_{ijk}$ over the edges of $T$. Because of the antisymmetry of the $y_{ijk}$'s, signs cannot be avoided here. In fact, the
correspondence between monomials and sub-$3$-graphs of $\Gamma_m$
is not one-to-one and a sub-$3$-graph
determines a monomial only up to sign. But these signs can be fixed unambiguously, although we won't explain this here (see \cite{MV1,MV2}).

If $m$ is even, then  one has $
\P_m=0
$, because there are no spanning trees in $3$-graphs with even number of vertices. If $m$ is odd, then $\P_m$ is a homogeneous polynomial of degree $(m-1)/2$ in the $y_{ijk}$'s. For example, one has $$\P_3=y_{123}$$ (the $3$-graph $\Gamma_3$ with three vertices and one edge is itself a tree). If $m=5$, we have
\begin{equation} \label{m=5}
\P_5=   y_{123}\, y_{145} -  y_{124}\, y_{135} +  y_{125}\, y_{134}
\  \pm \ \ldots~,
\end{equation}
where the right-hand side is a sum of $15$ similar terms corresponding
to the $15$ spanning trees of $\Gamma_5$.
If we visualize the edges of $\Gamma_m$ as  Y-shaped
objects \MILN, then the spanning tree corresponding to the first term 
of~(\ref{m=5})  
will look like on Figure~\ref{figexi}.  

We can now state one of the main results of \cite{MV1,MV2}.

 \begin{theorem}{\rm\cite{MV1,MV2}}\label{thm:alex-trees}\qua
Let $L$ be an algebraically split oriented link with
$m$ components. Then 
  \begin{equation}
  \label{eq:sqtrees}
  c_{2m-2}(L)= \bigl(\P_m(\mu_{ijk}(L))\bigr)^2~,
  \end{equation} where $\P_m(\mu_{ijk}(L))$ means the result of evaluating the polynomial $\P_m$ at $y_{ijk}=\mu_{ijk}(L)$.
\end{theorem}

For $m=3$, we find again Cochran's formula \eqref{Coc}, but for $m\ge 5$ our formula is new. For example, when $m=5$, 
we obtain that the first non-vanishing coefficient of $\nabla_L(z)$ 
for algebraically split links with $5$ components is equal to 
\begin{eqnarray*} 
&c_8(L)\hspace{-6pt}&= \P_5(\mu_{ijk}(L))^2\\
&&=
\bigl( \mu_{123}(L)\mu_{145}(L) - \mu_{124}(L)\mu_{135}(L)
+ \mu_{125}(L)\mu_{134}(L)  \ \pm \ \ldots\bigr)^2, 
\end{eqnarray*}
where $\P_5(\mu_{ijk}(L))$ consists of $15$ terms corresponding to
the spanning trees  of $\Gamma_5$.

\section{A proof via the Pfaffian Matrix-Tree Theorem of \cite{MV1}}

The first proof of Theorem~\protect{\ref{thm:alex-trees}} was given in \cite{MV1}. One of the main results of that paper is a \emph{Pfaffian Matrix-Tree
Theorem} which is the analog for $3$-graphs of the classical Matrix-Tree Theorem (see Section \ref{sec:2}). It expresses the generating function of spanning trees  on a $3$-graph $G$ as the Pfaffian of a matrix $\Lambda(G)^{(p)}$ associated to $G$. 

If $G$ is the complete $3$-graph $\Gamma_m$, this theorem says the following.

\begin{theorem}{\rm\cite{MV1}}\label{pmt}\qua
The generating function of spanning trees on the complete $3$-graph
$\Gamma_m$ is given by $$\P_m= (-1)^{p-1}
\Pf(\Lambda(\Gamma_m)^{(p)})~,$$ where $\Lambda(\Gamma_m)$ is the
$m\times m$ skew-symmetric matrix with entries
$\Lambda(\Gamma_m)_{ij}=\sum_{k} y_{ijk},$ and $\Pf$ denotes the
Pfaffian.
\end{theorem}

Recall that the Pfaffian of a skew-symmetric matrix $A$ is a polynomial in the
coefficients of $A$ such that $$(\Pf A)^2=\det A~.$$
The matrix $\Lambda$ defined in \eqref{eq:skew} is obtained from $\Lambda(\Gamma_m)$ by substituting the triple Milnor number $\mu_{ijk}(L)$ for the indeterminate $y_{ijk}$. Hence,  Theorem~\ref{pmt} implies Theorem~\ref{thm:alex-trees}, since we know from Formula \eqref{eq:lev} that
$$
c_{2m-2}(L)=\det\LP = (\Pf \LP)^2~.
$$
For a definition of the matrix $\Lambda(G)$ and a statement of the Pfaffian Matrix-Tree
Theorem in the case of general $3$-graphs $G$, as well as for a proof, see \cite{MV1}. 

\section{A proof via finite type invariants \cite{MV2} }

As explained in Section \ref{sec:3}, the appearance of spanning trees in Formula~\eqref{eq:ho_tree} for the coefficient $c_{m-1}$ is very natural from the point of view of finite type invariants. A similar approach also leads to a proof of Theorem~\protect{\ref{thm:alex-trees}} via finite type invariants. This argument naturally generalizes to higher Milnor numbers. Let us briefly describe this approach (see \cite{MV2} for details).

The connection between the Alexander-Conway polynomial and the Milnor
numbers is established by studying their weight systems and then using
the Kontsevich integral. In the dual language of the space of chord
diagrams, the Milnor numbers correspond to the tree diagrams
(see~\cite{HM}) and the Alexander-Conway polynomial can be described
in terms of certain trees and wheel diagrams (see~\cite{KSA}
and~\cite{Vai}). However, for first non-vanishing terms, only tree 
diagrams matter, as the following Vanishing Lemma shows.

\begin{proposition}[Vanishing Lemma \cite{MV2}] \label{2.1} 
Let $D$ be a  degree-$d$ diagram 
 on \hbox{$m\geq 2$} solid circles, such that $D$ has no tree
 components of degree $\leq n-1$. Let $W_d$ be the Alexander-Conway weight system. 
If $d\leq n(m-1)+1$,  then $W_d(D)=0$ unless $D$ has exactly $m-1$
 components, each of which is a tree of degree $\geq n$.
\end{proposition}

This result is the generalization of the fact, shown in Section~\ref{sec:3}, that the Alexander-Conway weight system $W_d$ for $m$-component links is always zero in degrees $d<m-1$. Indeed, this fact is the $n=1$ case of the Vanishing Lemma. However, the proof in the general case is more complicated. It uses properties of the
Alexander-Conway weight system 
from~\cite{FKV} which are based on the connection between $\nabla$ and
the Lie superalgebra  $gl(1|1)$. 

In view of the relationship between Milnor numbers and tree diagrams studied in \cite{HM}, the Vanishing Lemma implies in a rather straightforward way the following result, which was first proved by Traldi~\cite{Tr2}
and Levine~\cite{Le2} using quite different methods. 

\begin{proposition}{\rm\cite{Tr2,Le2}}\label{5.2}\qua
Let $L$ be an oriented link such that all Milnor invariants of 
$L$ of degree\footnote{Here, the \emph{degree} of a Milnor  invariant is 
its Vassiliev degree, which is one less than 
its length (the number of its indices).
For example, linking numbers have degree one, and
triple linking numbers have degree two.} $\le n-1$ vanish. 
Then  for the coefficients 
$c_i(L)$ of the Alexander-Conway
polynomial 
$\nabla_L(z)=\sum_{i\ge 0}c_i(L)z^i$ 
we have 
\begin{itemize} 
        \item[(i)] \ $c_i(L)=0$ for $i<n(m-1),$
\item[(ii)] \ $c_{n(m-1)}(L)$ is a homogeneous polynomial $F_m^{(n)}$ of degree $m-1$ in the Milnor numbers of $L$ of degree $n$.
\end{itemize}
\end{proposition}

Using the approach via Seifert surfaces, Levine~\cite{Le2} (see also Traldi~\cite{Tr1,Tr2}) gives a formula for the polynomial $F_m^{(n)}$ as a determinant in the degree $n$ Milnor numbers of $L$. For $n=1$ and $n=2$,  this formula specializes to Formulas~\eqref{eq:ho_det} and~\eqref{eq:lev}, respectively.

From the point of view of the Alexander-Conway weight system, however, one is lead to an expression for the polynomial $F_m^{(n)}$ in terms of the spanning tree polynomials $\D_m$ and $\P_m$. 
Indeed, as explained in Section \ref{sec:3}, for $n=1$ the polynomial $F_m^{(1)}$ is easily recognized to be the spanning tree polynomial $\D_m$ in the linking numbers $\ell_{ij}$. How does this generalize to higher $n$?

Consider for example the case $n=2$, that is, the case of algebraically split links. Proposition~\ref{5.2}(ii) tells us that $c_{2m-2}(L)$ is a homogeneous polynomial $F_m^{(2)}$ of degree $m-1$ in triple Milnor numbers $\mu_{ijk}(L)$.

\begin{theorem}{\rm\cite{MV2}}\qua
  The polynomial $F_m^{(2)}$ is equal to $\P_m^2$, the square of the Pfaffian-tree polynomial $\P_m$.
\end{theorem}

\begin{figure}[ht!]
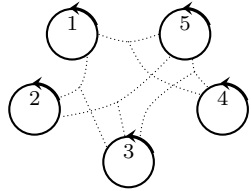
   
\begin{center}
\ExFive 

\caption{\label{FigExFive} 
A diagram contributing to $F_5^{(2)}$} 
\end{center}
\end{figure}

Here is a sketch of the proof. Triple Milnor numbers are dual to $\Y$-shaped diagrams \MILN, and the coefficients of the  polynomial $F_m^{(2)}$ can be computed from the Alexander-Conway weight system.  For example, the coefficient of the monomial $$\mu_{123}\,\mu_{145} \,\mu_{235}\,\mu_{345}$$ in $F_5^{(2)}$ is equal to the value of the Alexander-Conway weight system on the diagram in Figure~\ref{FigExFive}. 

\begin{figure}[ht!]
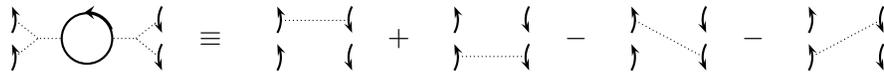

\begin{center}
$\Stepd  \ \ \ \equiv \    \Stepdii \ \ \ + \Stepdiia \ \ \ -\Stepdiib\ \ \ - \Stepdiic$
\end{center}\caption{\label{relAC} An identity modulo the Alexander-Conway relations}
\end{figure}

The coefficients of $F_m^{(2)}$  can be computed recursively by the relation in Figure~\ref{relAC}, which follows from identities proved in \cite{FKV}.
\vskip 8pt

Indeed, the relation in Figure~\ref{relAC} together with the smoothing relation~(\ref{smooth}) allows one  to reduce a diagram consisting of $m-1$ $\Y$'s on $m$ solid circles to a linear combination of diagrams consisting of $m-3$ $\Y$'s on $m-2$ solid circles. (We are leaving out some details here.) This gives recursion formulas expressing $F_m^{(2)}$ in terms of $F_{m-2}^{(2)}$. 

 Let us state an example of such a recursion formula. It is convenient to write the antisymmetric triple Milnor number formally as an exterior  product $$\mu_{ijk}=v_i\wedge v_j\wedge v_k$$ 
and to consider $F_m^{(2)}$ as an 
expression in the indeterminates  $v_i$: 
 $$ F_m^{(2)}=F_m^{(2)}(v_1,v_2,\ldots, v_m)~.$$ 
Then the relation in Figure~\ref{relAC} implies for example that $F_m^{(2)}$ satisfies the recursion relation
\begin{align*}
\left[\frac {\partial^2 \, F_m^{(2)}} {\partial\mu_{123}\,
\partial\mu_{145}}\right]_{v_1=0} &= F_{m-2}^{(2)}( v_3+v_4,v_2, v_5, \ldots) 
+ F_{m-2}^{(2)}(v_2 +v_5, v_3, v_4, \ldots)\\ 
& \ -F_{m-2}^{(2)}( v_2+v_4,v_3, v_5, \ldots) - 
 F_{m-2}^{(2)}( v_3 +v_5, v_2, v_4, \ldots) ~.
\end{align*}

\begin{figure}[ht!]
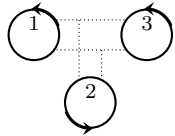
   
\begin{center}
\ExThree 
\end{center}\caption{\label{ext3} The only diagram contributing to
$F_3^{(2)}$} 
\end{figure}

It turns out that this and similar recursion relations  are enough to determine the polynomial $F_m^{(2)}$ for all $m$, once one knows it for $m=2$ and $m=3$. But it is easy to see that $F_{2}^{(2)}=0$, while $F_{3}^{(2)}=\mu_{123}^2$, the only non-zero diagram contributing to $F_{3}^{(2)}$ being the diagram in Figure~\ref{ext3}.

We now claim that this implies that $F_m^{(2)}$ is equal to $\P_m^2$, the square of the Pfaffian-tree polynomial $\P_m$. Clearly this is true for $m=2$ and $m=3$, and the proof consists of showing that $\P_m^2$ satisfies the same recursion relations as $F_m^{(2)}$. This uses the following two relations \eqref{eq:recursion2} and \eqref{eq3t} satisfied  by $\P_m$ itself, which follow more or less directly from the definition of $\P_m$ as the spanning tree generating function of the complete $3$-graph $\Gamma_m$ (see \cite{MV1}).

The first is a \emph{contraction-deletion} relation
\begin{equation}
  \label{eq:recursion2}
  \P_m= y_{123}\P_{m-2}(v_1+v_2+v_3,v_4,\ldots,v_m) +
  \left[\P_m\right]_{ y_{123}=0}~.
\end{equation}
Here, we have again written the indeterminate $y_{ijk}$ as an exterior product $v_i\wedge v_j\wedge v_k~$. The first term on the right hand side corresponds to the spanning trees on $\Gamma_m$ containing the edge $\{1,2,3\}$, and the second term to those that do not.  Note that a similar contraction-deletion
 relation exists for the classical spanning tree generating function  for usual graphs. 

The second relation is called {\em Three-term relation} in \cite{MV1}. It states that 
 \begin{equation}  \label{eq3t}
 \P_m(v_2+v_3,v_4,\ldots) +\P_m(v_3+v_4,v_2,\ldots)
 +\P_m(v_2+v_4,v_3,\ldots) =0 
 \end{equation}
  where the dots stand for $v_5,v_6,\ldots, v_{m+2}$.

The contraction-deletion relation and the three-term relation imply, by some algebraic manipulation, that $\P_m^2$ satisfies the same recursion relations as $F_m^{(2)}$.

Thus, although the recognition of the polynomial  $F_m^{(2)}$ as being equal to  the squared spanning tree polynomial $\P_m^2$ is not quite as immediate from the Alexander-Conway weight system as the identification of $F_m^{(1)}$ with the spanning tree polynomial $\D_m$ in Section \ref{sec:3}, it is still quite natural. Indeed, it is based on the fact that the recursion relations have two natural interpretations, one coming from the weight system relations in Figure~\ref{relAC}, and one coming from the contraction-deletion relation and the three-term relation for the Pfaffian-tree polynomial $\P_m$.

% This concludes our discussion of the proof of Theorem~\ref{thm:alex-trees} 
% via finite type invariants. 

\section{Some generalizations to higher Milnor numbers} 
The polynomial $F_m^{(n)}$ can be determined explicitly for higher values of $n$ also. The answer can be expressed in terms of the spanning tree polynonmials $\D_m$ or $\P_m$. One obtains the following result for links with vanishing Milnor numbers up to a given degree.

\begin{theorem}{\rm\cite{MV2}}\label{thm:general}\qua
Let $L$ be an oriented $m$-component link 
with vanishing Milnor numbers of degree $p <n$
and let 
$\nabla_L(z)=\sum_{i\geq 0} c_i(L) z^i $
be its Alexander-Conway polynomial. Then
$c_i=0$ for $i<n(m-1)$ 
and
\begin{equation*}\label{eq:genform}
c_{n(m-1)}(L)= 
\begin{cases} \ \ \D_m(\ell^{(n)}_{ij})~,   &\text{if $n$ is odd}\\                       (\P_m(\mu^{(n)}_{ijk}))^2\, , &\text{if $n$ is even},                 \end{cases}
\end{equation*}
where $\ell^{(n)}_{ij}$ and $\mu^{(n)}_{ijk}$ are certain linear
combinations of the Milnor numbers of $L$ of degree $n$.
\end{theorem}

\begin{figure}[ht!]
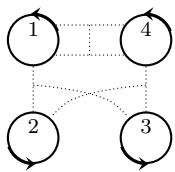
   
\begin{center}
\ExFouri

\caption{\label{figFour} A diagram contributing to $G_4^{(2)}$}
\end{center}
\end{figure}

Note that if $m$ is even then $\P_m=0$ and so if $n$ is also even, then 
the coefficient $c_{n(m-1)}(L)$ is always zero. In this case the 
Vanishing Lemma~\ref{2.1} leads to an expression for the next coefficient $c_{n(m-1)+1}(L)$  in terms of a certain polynomial $G_m^{(n)}$. (See Figure~\ref{figFour} for an example of a diagram contributing to $G_4^{(2)}$.) 

This polynomial can again be expressed via spanning trees (see \cite{MV2}).

\Addresses\recd

\end{document}